\documentclass[11pt]{amsart}
\usepackage{amscd}
\usepackage{amsfonts}
\usepackage{color}
\usepackage{amsmath,amssymb,amsthm,latexsym}
\usepackage{amscd}
\usepackage{hyperref}
\usepackage{enumerate}

\newtheorem{theorem}{Theorem}[section]
\newtheorem{proposition}[theorem]{Proposition}
\newtheorem{definition}[theorem]{Definition}
\newtheorem{corollary}[theorem]{Corollary}
\newtheorem{lemma}[theorem]{Lemma}
\usepackage{amsmath}
\usepackage{amsfonts}
\usepackage{amsmath,amsthm}
\usepackage{amscd}
\usepackage[all]{xy}
\numberwithin{equation}{section}
\theoremstyle{remark}
\newtheorem{remark}[theorem]{Remark}
\newtheorem{example}[theorem]{\bf Example}
\newcommand{\R}{\mathbb{R}}
\newcommand{\C}{\mathbb{C}}
\newcommand{\D}{\mathbb{D}}
\newcommand{\ee}{\mathrm{e}}
\newcommand{\dd}{\mathrm{d}}

\newcommand{\FC}{\mathcal{C}}

\def\<{{\langle}}
\def\>{{\rangle}}
\begin{document}

\title[Willmore surfaces via loop groups: a survey]{\bf{Willmore surfaces in spheres via loop groups: a survey}}
\author{Josef F. Dorfmeister, Peng Wang}

\date{\today}
\maketitle

\begin{abstract}In the past decades, the authors made some systematic research on global and local properties of Willmore surfaces in terms of the DPW method. In this note we
give a survey, mainly including the basic framework of the DPW method for the global geometry of Willmore surfaces via the conformal Gauss map, applications on constructions of Willmore $2$-spheres, characterizations of minimal surfaces, Willmore deformations of Willmore surfaces and Bjoerling problems for Willmore surfaces. Moreover, we also obtained some results on harmonic maps via DPW, including a duality theorem for harmonic maps into an inner non-compact symmetric space and its dual inner compact symmetric space, and harmonic maps of finite uniton type.
\end{abstract}

\vspace{0.5mm}  {\bf \ \ ~~Keywords:}
Willmore surface; conformal Gauss map; normalized potential; non-compact symmetric space; Iwasawa decomposition.   \vspace{2mm}

{\bf\   ~~ MSC(2010): \hspace{2mm} 53A30, 53C30, 53C35}

\tableofcontents

\section{Introduction}
Willmore surfaces play an important role in the conformal geometry of surfaces in space forms. The study can be traced back to the French mathematician Germain's work on plate theory in 1820s. In the 1920s, Blaschke and his school began a systematic study of Willmore surfaces (where they called them conformal minimal surfaces), which was written mainly in German and attracted not too much attention. In the 1960s, Willmore proposed the famous Willmore conjecture for 2-tori, which turned out to be an important open problem in global differential geometry. The study of the Willmore conjecture led to much progress in many aspects of geometry, for example,
by the famous work of Li-Yau \cite{Li-Yau} on conformal area and the
celebrated work of Marques and Neves on geometric measure theory \cite{Marques}.  In particular, Bryant made several important contributions to the understanding of the geometry of
Willmore surfaces in the 3-sphere \cite{Bryant1984}. Later his work was generalized to Willmore surfaces in an n-sphere by Ejiri \cite{Ejiri1988}. In particular, it was shown that a surface is Willmore if and only if its conformal Gauss map is a harmonic map into some Grassmannian. It is well-known that such harmonic maps admit a loop group representation, that is, an application of the DPW method. It is therefore a  natural idea to deal with Willmore surfaces via the DPW method (see for example \cite{Helein}, \cite{Xia-Shen}). In a series of papers we re-visited this idea and made some progress on the global treatment of Willmore surfaces via the DPW method.

Although deriving the conformal Gauss map from a Willmore surface is straightforward, the converse way is not so. In fact it was very surprising to these authors that, to the best of our knowledge, there did not exist in the literature any characterization of those harmonic maps which can serve as the conformal Gauss map of a Willmore surface. In \cite{DoWa11}, we
provided such a theorem. With this theorem at hand, we can naturally apply the DPW method to Willmore surfaces,  and some simple surfaces can be constructed easily in this way \cite{Wang-1}. Also, several classes of Willmore surfaces can be discussed and constructed fairly explicitly
following our approach. A particularly interesting class of such examples consists of  the Willmore 2-spheres  in $S^{n+ 2}.$  For such a class of examples we need to generalize the work of Uhlenbeck
on harmonic two-spheres in $U(n)$ and of Burstall-Guest on harmonic two-spheres in general compact, inner symmetric spaces $G/K$.

When carrying out this program, we realized that there is a natural duality relating  harmonic maps into a non-compact inner symmetric space  to harmonic maps into its compact dual inner symmetric space. Locally the converse also holds. Applying the duality theorem permits us to translate the main result of Burstall-Guest for harmonic 2-spheres in compact inner symmetric spaces to a result about  harmonic 2-spheres in non-compact inner symmetric spaces.
This yields a coarse classification of Willmore 2-spheres in $S^{n+2}$ and leads to constructions of many new Willmore 2-spheres.

Independent  of the characterization  of Willmore surfaces of a
specific type, also Willmore surfaces with symmetries are of great interest. These authors have investigated the various types of symmetries and the possibilities to obtain Willmore surfaces with specific symmetries (both for orientation preserving resp. reversing symmetries) in \cite{DoWa-sym1} and \cite{DoWa-sym2}.
To keep the paper reasonably short, in this note we will only present basic principles and some examples.

The note is organized as follows:  We will recall the basic surface theory of Willmore surfaces, following the treatment of Burstall-Pedit-Pinkall  \cite{BPP} in Section 2. Then we will introduce the characterization of Willmore surfaces via harmonic maps  in Section 3, and combine these results with the DPW method in Section 4. We will see how one can describe minimal surfaces in space forms under this treatment in Section 5. Section 6 is devoted to  classifications and constructions of a special class of Willmore 2-spheres in $S^6$. Then in Section 7 we will consider symmetric Willmore surfaces. Section 8 concerns the Bjoerling problem for Willmore surface. In Section 9, we will collect some progress on the harmonic map theory which we developed in our study of Willmore surfaces. Finally we end the note with some open problems in Section 10.

\section{Willmore surfaces in spheres}

In this section we recall our formulation \cite{DoWa11} of the treatment \cite{BPP} of Willmore surfaces
in $S^{n+2}$ via the projective light cone model of the conformal geometry of $S^{n+2}$.

\subsection{ Basic setting for a description of harmonic maps into spheres}  Let $\mathbb{R}^{n+4}_1$ denote the Minkowski space, i.e. $\mathbb{R}^{n+4}$ equipped with the Lorentzian metric
$$\langle x,y\rangle=-x_{0}y_0+\sum_{j=1}^{n+3}x_jy_j=x^t I_{1,n+3} y,\ \hspace{3mm}  I_{1,n+3}=diag(-1,1,\cdots,1).$$
The projective light cone $
Q^{n+2}=\{\ [x]\in\mathbb{R}P^{n+3}\ |\   x \in \mathcal{C}_+^{n+3} \}$
with the induced conformal metric, is conformally equivalent to $S^{n+2}$. Here $\mathcal{C}^{n+3}_+=\lbrace x \in \mathbb{R}^{n+4}_{1} |\<x,x\>=0 , x_0>0 \}.$
Moreover, the oriented conformal group of
 $S^{n+2}$ is exactly $SO^+(1,n+3)$, with $SO^+(1,n+3)$ being the connected component containing identity of the special orthogonal group
of $\R^{n+4}_{1}$. The group acts on $S^{n+2}\cong Q^{n+2}$ by
$T([y])=[Ty],\ T\in SO^+(1, n+3)$.

Let $y:M\rightarrow S^{n+2}$ be a conformal immersion from a Riemann surface $M$, and  $U\subset M$  a contractible open subset. A local canonical lift of $y$ is a map $Y:U\rightarrow \mathcal{C}_+^{n+3} $ such that $\pi\circ Y=y$ and $|{\rm d}Y|^2=|{\rm d}z|^2$, with  $z$ a local complex coordinate on $U$.
We consider the trivial global vector bundle $M\times \mathbb{R}^{n+4}_{1}$ over $M$ and restrict it to $U$. We put
\begin{equation}
{V}|_{p}:={\rm Span}_{\mathbb{R}}\{Y,{\rm Re}Y_{z},{\rm Im}Y_{z},Y_{z\bar{z}}\}|_{z=z(p)}
\end{equation}
for $p\in U$. It is easy to see that $V_p$ is independent of the choice of $z$ and it is  an oriented rank-4 Lorentzian subspace of $\R^{n+4}_1$.
Moreover, if one chooses a covering of $M$ by $U$'s  as above, then the above defined bundles  over the $U$'s combine to yield a global bundle over $M$. Furthermore, the bundle decomposition $M\times \mathbb{R}^{n+4}_{1}=V\oplus V^{\perp}$ with $V^{\perp}|_p \perp V_p$ is globally conformally invariant. The complexifications of  $V$ and $V^{\perp}$ will be  denoted by $V_{\mathbb{C}}$ and $V^{\perp}_{\mathbb{C}}$ respectively. There exists a unique $N\in\Gamma(V|_U)$ on $U$     such that
\[ \langle N,Y_{z}\rangle=\langle N,Y_{\bar{z}}\rangle=\langle
N,N\rangle=0,\langle N,Y\rangle=-1.\]
Note that $N\equiv 2Y_{z\bar{z}}\!\!\mod Y$ and  ${V}|_{p}:={\rm Span}_{\mathbb{R}}\{Y,{\rm Re}Y_{z},{\rm Im}Y_{z},N \}|_{z=z(p)}.$
We define \emph{the conformal Gauss map} of $y$.
\begin{definition} \label{def-gauss} $($\cite{Bryant1984},\cite{BPP},\cite{Ejiri1988}$)$
Let $y:M\to S^{n+2}$ be a conformal immersion from a Riemann surface $M$. The \emph{conformal Gauss map} of $y$ is defined by
\begin{equation}\begin{array}{ccccc}
                { Gr_y :} & M &\rightarrow&
Gr_{1,3}(\mathbb{R}^{n+4}_{1}) &= SO^+(1,n+3)/SO^+(1,3)\times SO(n)\\
                \ & p\in M & \mapsto & V_p &\ \\
                \end{array}
\end{equation}
Locally we can express $Gr_y$  as follows \cite{Ma}:
$
Gr_y=Y\wedge Y_{u}\wedge Y_{v}\wedge N=-2i\cdot Y\wedge Y_{z}\wedge
Y_{\bar{z}} \wedge N$, with $z = u + i v$. For simplicity we will frequently write just $Gr$ instead of $Gr_y.$
\end{definition}

\subsection{Structure equations and integrability conditions for conformal maps}
The structure equations and the integrable equations of $y$ are (\cite{BPP}, \cite{Ma2006}):
\begin{equation}\label{eq-moving}
\left\{\begin {array}{lllll}
Y_{zz}=-\frac{s}{2}Y+\kappa,\\
Y_{z\bar{z}}=-\langle \kappa,\bar\kappa\rangle Y+\frac{1}{2}N,\\
N_{z}=-2\langle \kappa,\bar\kappa\rangle Y_{z}-sY_{\bar{z}}+2D_{\bar{z}}\kappa,\\
\psi_{z}=D_{z}\psi+2\langle \psi,D_{\bar{z}}\kappa\rangle Y-2\langle
\psi,\kappa\rangle Y_{\bar{z}},
\end {array}\right.
\end{equation}
and
\begin{equation}\label{eq-integ}
\left\{\begin {array}{lllll} \frac{1}{2}s_{\bar{z}}=3\langle
\kappa,D_z\bar\kappa\rangle +\langle D_z\kappa,\bar\kappa\rangle, & \hbox{Gauss equation}\\
{\rm Im}(D_{\bar{z}}D_{\bar{z}}\kappa+\frac{\bar{s}}{2}\kappa)=0,& \hbox{Codazzi equation }\\
R^{D}_{\bar{z}z}=D_{\bar{z}}D_{z}\psi-D_{z}D_{\bar{z}}\psi =
2\langle \psi,\kappa\rangle\bar{\kappa}- 2\langle
\psi,\bar{\kappa}\rangle\kappa, & \hbox{Ricci equations}
\end {array}\right.
\end{equation}
respectively. Here  $D$ {denotes}  the normal connection and $\psi\in\Gamma(V^{\perp})$. The section $\kappa\in \Gamma(V_{\mathbb{C}}^{\perp})$ is called \emph{the conformal Hopf differential} of $y$,
and the function $s$ is called \emph{the Schwarzian} of $y$.
We refer to \cite{BPP} for a more detailed discussion of these complex valued functions/maps.

Let $\{\psi_j,j=1,\cdots,n\}$ be an oriented orthonormal frame  of $V^{\perp}$  with $D_z\psi_j=\sum_{1\leq l\leq n}b_{jl}\psi_l$ and $b_{jl}+b_{lj}=0.$
Set  $\kappa=\sum_{j=1}^{n}k_j\psi_j$, $k=|\kappa|=\sqrt{\sum_{j=1}^{n}|k_j|^2}$ and $D_{\bar{z}}\kappa=\sum_{1\leq j\leq n}\beta_j\psi_j$ with $\beta_j=k_{j\bar{z}}- \sum_{1\leq j\leq n}\bar{b}_{jl}k_l$. Set
$\phi_1=\frac{1}{\sqrt{2}}(Y+N),\ \phi_2=\frac{1}{\sqrt{2}}(-Y+N),\ \phi_3= Y_z+Y_{\bar{z}},\ \phi_4=i(Y_z-Y_{\bar{z}})$. Then  $
F:=\left(\phi_1,\phi_2,\phi_3,\phi_4,\psi_1,\cdots,\psi_n\right)$ is  a frame of  $Gr$ and we have the following result.
\begin{proposition}\label{frame}\
\begin{enumerate}
    \item
The Maurer-Cartan form $\alpha=F^{-1}\dd F$ of the frame $F$ is of the form \[\alpha=F^{-1}\dd F=\left(
                   \begin{array}{cc}
                     A_1 & B_1 \\
                     B_2 & A_2 \\
                   \end{array}
                 \right)\dd z+\left(
                   \begin{array}{cc}
                     \bar{A}_1 & \bar{B}_1 \\
                     \bar{B}_2 & \bar{A}_2 \\
                   \end{array}
                 \right)\dd\bar{z},\]
with
\begin{equation}A_1=\left(
                             \begin{array}{cccc}
                               0 & 0 & s_1 & s_2\\
                               0 & 0 & s_3 & s_4 \\
                               s_1 & -s_3 & 0 & 0 \\
                               s_2 & -s_4 & 0 & 0 \\
                             \end{array}
                           \right),\   A_2=\left(
                             \begin{array}{cccc}
                               b_{11} & \cdots &  b_{n1} \\
                               \vdots& \vdots & \vdots \\
                               b_{1n} &\cdots & b_{nn} \\
                             \end{array}
                           \right),\ \end{equation}
\begin{equation}\label{s}
 s_1=\frac{1-s-2k^2}{2\sqrt{2}},\ s_2=-\frac{i(1+s-2k^2)}{2\sqrt{2}},\ s_3=\frac{1+s+2k^2}{2\sqrt{2}},\ s_4=-\frac{i(1-s+2k^2)}{2\sqrt{2}},
\end{equation}
and
\begin{equation} \label{B1}
B_1=\left(
      \begin{array}{ccc}
         \sqrt{2} \beta_1 & \cdots & \sqrt{2}\beta_n \\
         -\sqrt{2} \beta_1 & \cdots & -\sqrt{2}\beta_n \\
        -k_1 & \cdots & -k_n \\
        -ik_1 & \cdots & -ik_n \\
      \end{array}
    \right),  \ \
B_2=-B_1^tI_{1,3}=\left(
      \begin{array}{cccc}
        \sqrt{2} \beta_1& \sqrt{2} \beta_1 & k_1& i k_1 \\
        \vdots & \vdots & \vdots & \vdots\\
      \sqrt{2} \beta_n& \sqrt{2} \beta_n & k_n & ik_n \\
      \end{array}
    \right).
 \end{equation}
 \item
 Conversely, assume the frame $F=(\phi_1,\cdots,\phi_4,\psi_1,\cdots,\psi_{n+4}):U\rightarrow SO^+(1,n+3)$ has its Maurer-Cartan form $\alpha=F^{-1}\dd F$ of the form as above. Then
\begin{equation}\label{eq-y-from-map}y=\pi_0(F)=:\left[ (\phi_1-\phi_2)\right]
\end{equation}
is a conformal immersion from $U$ into $Q^{n+2}\cong S^{n+2}$ with canonical lift $\frac{1}{\sqrt{2}}(\phi_1-\phi_2)$.
\end{enumerate}
    \end{proposition}

\begin{corollary}
The equation $B_1 ^t I_{1,3} B_1 = 0$ holds for $Y$. In particular, $Rank(B_1)\leq 2$.
\end{corollary}


\vspace{2mm}
The conformal Gauss map $Gr$ induces {on $M$ a conformally invariant (possibly degenerate) metric  (see e.g. \cite{encyc} for a definition))
$$
g:=\frac{1}{4}\langle {\rm d}G,{\rm d}G\rangle=\langle
\kappa,\bar{\kappa}\rangle|\dd z|^{2}$$
 which degenerates exactly at umbilical points of $y$, i.e.. the points where  $\kappa$ vanishes.  (see \cite{BPP}).

\subsection{{Definition and basic results ablout Willmore surfaces in spheres}}
{Next we define the Willmore functional.}
\begin{definition}$($\cite{BPP}, \cite{Ma2006}$)$ The \emph{ Willmore functional} of $y$ is
defined as follows via the metric defined just above:
\begin{equation}\label{eq-W-energy}
W(y):=2i\int_{M}\langle \kappa,\bar{\kappa}\rangle dz\wedge
d\bar{z}.
\end{equation}
The surface $y$ is called a
\emph{Willmore surface}, if it is a critical surface w.r.t the Willmore functional. Note that this definition of the Willmore functional coincides with the usual definition.
\end{definition}
 Moreover, Willmore surfaces can be characterized as follows
\cite{Blaschke}, \cite{Bryant1984}, \cite{BPP}, \cite{Ejiri1988}, \cite{Rigoli1987}, \cite{Wang1998}.
\begin{theorem}\label{thm-willmore} The following three conditions
are equivalent:
 \begin{enumerate}
\item $y$ is Willmore;
\item The conformal Gauss map $Gr$ is a conformally harmonic map into
$G_{3,1}(\mathbb{R}^{n+3}_{1})$;
\item The conformal Hopf differential $\kappa$ of $y$ satisfies the Willmore equation:
\begin{equation}\label{eq-willmore}
D_{\bar{z}}D_{\bar{z}}\kappa+\frac{\bar{s}}{2}\kappa=0.
\end{equation}
 \end{enumerate}
\end{theorem}
 The loop group method outlined below
permits to construct all Willmore surfaces, at least in principle. For Willmore surfaces and  minimal surfaces in a space form we have the following well known results.
\begin{example}
 \begin{enumerate}
\item It is well-known that all minimal surfaces in Riemannian space forms  can be considered to be  Willmore surfaces in some $S^{n+2}$ (\cite{Bryant1984}, \cite{Weiner}).  In particular, let  $M$ be a compact Riemann surface and $x:M\setminus\{p_1,\cdots, p_r\}\rightarrow \mathbb{R}^{n+2}$ be a complete minimal surface with  embedded planar ends $\{p_1,\cdots, p_r\}$. By the inverse of the stereographic projection $\pi$, we get a Willmore immersion $y=\pi^{-1}(x): M\rightarrow S^{n+2}$ (\cite{Bryant1984}, \cite{Bryant1988}).

\item The isotropic (or super-conformal \cite{Bryant1982}) surfaces in $S^4$ provide another important kind of Willmore surfaces. It is well known (\cite{Ejiri1988}, \cite{Mus1}, \cite{Mon}, \cite{BFLPP}) that
 a Willmore  immersion $y:S^2\rightarrow S^4$ is either an istropic $2$-sphere or coming from a complete minimal surface of genus 0 in $\mathbb{R}^4$ with embedded planar ends.
\item The first non-minimal Willmore surface was constructed by Ejiri in \cite{Ejiri1982}, which is a homogeneous torus in $S^5$. Using the Hopf bundle, Pinkall produced a family of non-minimal Willmore tori in $S^3$ via elastic curves (\cite{Pinkall1985}).
 \end{enumerate}
\end{example}
\subsection{{ Transforms of Willmore surfaces and $S-$Willmore surfaces}}

 Transforms  of Willmore surfaces lead to  new Willmore surfaces and to a deeper  understanding of Willmore surfaces. The most simple  transform is the construction of a  ``dual surface".

\begin{definition}\label{def-dual} $($\cite{Bryant1984}, page 399 of \cite{Ejiri1988}$)$
Let $y:M\rightarrow S^{n+2}$ be a Willmore surface with $M_0$ being the set of umbilical points of $y$. The conformal map $\hat{y}:M\setminus M_0\rightarrow S^{n+2}$ is called a ``dual surface" of $y$, if either $\hat y$ reduces to a point, or $\hat y$  is an immersion  on an open dense subset and $\hat{y}$ has the same conformal Gauss as  $y$ {but} with opposite orientation (Hence $\hat{y}$   is also a Willmore surface).
\end{definition}
In general a Willmore surface in $S^{n+2}$ may not admit a dual surface. To describe Willmore surfaces with dual surfaces, Ejiri introduced the so-called {\em S-Willmore surfaces}  in
 \cite{Ejiri1988}.
 It is convenient to define S-Willmore surfaces as follows (see also \cite{Ma2005}):
\begin{definition} $($\cite{Ejiri1988}$)$
A Willmore immersion $y:M\rightarrow S^{n+2}$ is called an S-Willmore surface if on any  contractible open subset $U$ away from the umbilical points, the conformal Hopf differential $\kappa$ of $y$ satisfies
$D_{\bar{z}}\kappa || \kappa,\
~\hbox{ i.e.  }\  D_{\bar{z}}\kappa+\frac{\bar{\mu}}{2}\kappa=0~\hbox{ for some } \mu:U\rightarrow \mathbb{C}. $
\end{definition}
By definition and Proposition \ref{frame}, we have the following observation.
\begin{corollary}\label{S-frame} Let $y$ be a non totally umbilical Willmore surface. Then $y$ is S-Willmore  if and only if the (maximal) rank of $B_1$ in Proposition \ref{frame} is $1$.
\end{corollary}

\begin{theorem}\label{thm-dual gauss map} \
 \begin{enumerate}
 \item $($\cite{Bryant1984}, Theorem 7.1 of \cite{Ejiri1988}$)$ A non totally umbilical Willmore surface $y$ is S-Willmore  if and only if it has a unique dual  surface $\hat y$ on $M\setminus M_0$. Moreover, if $y$ is S-Willmore, the dual map $\hat y$ can be extended to $M$.
\item $($Theorem 2.9 of \cite{Ma}$)$ If the dual surface $\hat y$ of $y$ is immersed at $p\in M$,  then $Gr_{\hat{y}}(p)$ spans the same subspace as $Gr_y(p)$, but its orientation  is  opposite to the one of $Gr_y$.
\item
 \cite{Ejiri1988}, \cite{Ma2005} \label{Conformal-Gauss-map} Let $f:M\rightarrow SO^+(1,n+3)/SO^+(1,3)\times SO(n)$ be  the oriented   conformal Gauss map of a non totally umbilical Willmore surface $y:M\rightarrow S^{n+2}$ such that $f$ is also the conformal Gauss map of a surface $\tilde{y}$ on an open subset $U$ of $M$. Then
$\tilde{y}=y$ on $U$.

\item
$($\cite{Helein}, \cite{Ma-W1}$)$ A Willmore surface $y$ is conformally equivalent to a minimal surface in $R^{n+2}$  if and
only if its conformal Gauss map $Gr$ contains a constant lightlike vector. Here we say ``the conformal Gauss map contains a constant lightlike vector $Y_0 $'' if there exists a non-zero constant lightlike vector $Y_0$ in $\mathbb{R}^{n+4}_1$ satisfying $Y_0\in V_{p}$ for all $p\in M$.
\end{enumerate}
\end{theorem}

\subsection{{About Willmore surfaces with branch points}}

There exist Willmore surfaces which fail to be immersions at some points or lines. To include surfaces of this type, we introduce the notion of {\em Willmore maps} and {\em strong Willmore maps}.

\begin{definition}
A smooth map $y$ from a Riemann surface $M$ to $S^{n+2}$ is called a Willmore map if it is a conformal Willmore immersion on an open dense subset $\hat{M}$ of $M$.  If $\hat{M}$ is maximal, then the points in $M_0=M\backslash \hat{M}$ are called branch points of $y$, at which points $y$ fails to be an immersion. Moreover, $y$ is called a strong Willmore map if it is a Willmore map and if the conformal Gauss map $Gr: \hat{M}\rightarrow SO^+(1,n+3)/SO^+(1,3)\times SO(n)$ of $y$ can be extended smoothly (and hence real analytically) to $M$.
\end{definition}

\begin{definition} For a conformal map $f:M\rightarrow SO^+(1,n+3)/SO^+(1,3)\times SO(n)$,
We define the linear map $\mathfrak{B}$ as follows:
\begin{equation}\label{eq-B}
 \begin{array}{cccc}
  \mathfrak{B}:&\R^{n+4}_1& \rightarrow &\R^{n+4}_1\\
 \ & X &\mapsto& F \alpha_{\mathfrak{p}}'(\frac{\partial}{\partial z})F^{-1}X
 \end{array}
\end{equation}
It is straightforward to see that $\mathfrak{B}$ is well-defined, i.e., independent of the choice of the  frame  $F$.
\end{definition}
\begin{remark}\
\begin{enumerate}
\item Note that $\mathfrak{B}(V_{\C})=\partial f$, if we use the notion $\partial f$ as Chern and Wolfson defined in \cite{Chern-W}.
\item In \cite{Hitchin}, Hitchin  introduced the idea of decomposing $\dd$ into the connection part and a nilpotent linear map. Here $\mathfrak{B}$ is the linear map. Note that the same idea is also used in the definition of $\beta$ in \cite{BR} if one uses $\beta$ acting on a section as in Proposition 1.1 of \cite{BR}.
\end{enumerate}
\end{remark}
Using the fact that
$\alpha_{\mathfrak{p}}'(\frac{\partial}{\partial z})=\left(
                   \begin{array}{cc}
                     0 & B_1 \\
                     -B_1^tI_{1,3} & 0 \\
                   \end{array}
                 \right),$
                 we obtain immediately
\begin{corollary}\
 \begin{enumerate}
\item $\mathfrak{B}(V_{\C})\subset V_{\C}^{\perp}$, $\mathfrak{B}(V_{\C}^{\perp})\subset V_{\C}$.
\item The linear map $\mathfrak{B}$ satisfies the ``restricted nilpotency condition"
\begin{equation}\label{eq-nilpotent-B}
\mathfrak{B}^2|_{V^\perp}=0.
\end{equation}
if and only if the Maurer-Carten form of $f$ satisfies
\[B_1^t I_{1,3}B_1=0.\]
\end{enumerate}
\end{corollary}

\section{A Weierstrass-Kenmotsu type theorem}
We define the strongly conformally harmonic map as follows.

\begin{definition}\cite{DoWa11} \label{stronglyconfharm}
Let $f: M\rightarrow SO^+(1,n+3)/SO^+(1,3)\times SO(n)$ be a harmonic map. We call $f$ a {\bf strongly conformally harmonic map} if  it satisfies the equation \eqref{eq-nilpotent-B}:  $\mathfrak{B}^2|_{V^\perp}=0.$
 \end{definition}
Then we have the following Weierstrass-Kenmotsu type theorem for $f$.
\begin{theorem} \cite{DoWa11}
Let $f: M\rightarrow SO^+(1,n+3)/SO^+(1,3)\times SO(n)$ be a non-constant strongly  conformally  harmonic map. Then either it contains a  constant lightlike vector, or up to a change of orientation, it is the oriented conformal Gauss map of a Willmore map $y$ on an open dense subset of $M$ where $y$ is immersed.
\end{theorem}

\begin{remark}
\begin{enumerate}
    \item We refer to \cite{DoWa11} for a more detailed discussion. Note that the key idea of the proof is that we can gauge the frame $F$ of the strongly conformally  harmonic map $f$ such that the Maurer-Cartan form has the desired form as in Proposition \ref{frame}. The harmonicity of $f$ makes it real analytic so that we can get the map $y=[Y]$ globally.
\item
    Note that some partial version of this theorem has been proved in \cite{Helein}, \cite{Xia-Shen}. The global version, together with a detailed discussion on the case of  $f$ containing a constant lightlike vector, is due to \cite{DoWa11}.
    \item
    In \cite{Bu19}, Burstall generalized our results to a local description of the conformal Gauss map of a conformal immersion in a bundle version. We refer to \cite{Bu19} for more details.
\end{enumerate}

\end{remark}

\section{The loop group method for Willmore surfaces}

In this section we recall briefly the loop group theory, the DPW method,  for harmonic maps and its application to  strongly conformally harmonic maps (\cite{DPW}, \cite{Wu}, \cite{Ba-Do}, \cite{DoWa12}).

\subsection{The loop group method for harmonic maps}
\subsubsection{Loop groups and decomposition theorems}
Let $G$ be a connected real matrix Lie group and $G^\C$ its complexification (\cite{Hochschild}). {Let $\ee$ denote the identity element of any of the groups occurring.} Let $\sigma$ be an inner involution of $G$ and $K$ a closed subgroup with $(Fix^\sigma(G))^0 \subset K \subset Fix^\sigma(G)$. Here the superscript ``${\ }^0$'' {denotes} the connected component containing the identity element of the topologic space.
Then $\sigma$ fixes $\mathfrak{k} = Lie K$.
The extension of $\sigma$  to an involution of $G^\C$ has
 $\mathfrak{k}^\C$  as  its fixed point algebra. Where needed, we define $K^\C$ to be the smallest complex subgroup of $G^\C$ which contains $K$.
 From harmonic maps into $G/K$, where $K^\C$ is the smallest complex subgroup of $G^\C$ which contains some $K$, we can derive all other harmonic map
 of type $G/K$.
It is known that $K^\C$ is in general not connected, but $Lie K^\C = (Lie K)^\C = \mathfrak{k}$ holds.

We recall the basic definitions about loop groups related to $G/K$.
{To define a Banach Lie group structure  on these groups we consider a weighted Wiener algebra $\mathfrak{A}_{\omega}$ on $S^1$, where $\omega$ denotes a ``symmetric weight of non-analytic type". Such a norm induces naturally a Banach matrix norm on the space of loops $Mat(m,\mathfrak{A}_{\omega})$, if
$G^\C \subset Mat(m,\C),$ and thus on the subgroup $\Lambda G^C$ of loops taking values in $G^\C$. Then the twisting homomorphism $\hat{\sigma}$ of $Mat(m,\mathfrak{A}_{\omega}),$ given by
$(\hat{\sigma}(\gamma))(\lambda) =
\sigma ( \gamma(-\lambda) )$ is continuous and the subgroups defined below inherit naturally the structure of a Banach Lie group. For more details see,
e.g. \cite{GoWa}, or mutatis mutandis the appendix of the paper
\cite{DoMaB2}.

 Then we put:
\begin{equation*}
\begin{array}{llll}
\Lambda G^{\mathbb{C}}_{\sigma} ~&=\{\gamma:S^1\rightarrow G^{\mathbb{C}}~|~ ,\
\sigma \gamma(\lambda)=\gamma(-\lambda),\lambda\in S^1  \},\\[1mm]
\Lambda G_{\sigma} ~&=\{\gamma\in \Lambda G^{\mathbb{C}}_{\sigma}
|~ \gamma(\lambda)\in G, \hbox{for all}\ \lambda\in S^1 \},\\[1mm]
\Omega G_{\sigma} ~&=\{\gamma\in \Lambda G_{\sigma}|~ \gamma(1)=\ee \},\\[1mm]
 \Lambda^{-} G^{\mathbb{C}}_{\sigma}  ~&=
\{\gamma\in \Lambda G^{\mathbb{C}}_{\sigma}~
|~ \gamma \hbox{ extends holomorphically to } |\lambda|>1 \cup\{\infty\} \},\\[1mm]
\Lambda_{*}^{-} G^{\mathbb{C}}_{\sigma} ~&=\{\gamma\in \Lambda G^{\mathbb{C}}_{\sigma}~
|~ \gamma \hbox{ extends holomorphically to } |\lambda|>1\cup\{\infty\},\  \gamma(\infty)=\ee \},\\[1mm]
\Lambda^{+} G^{\mathbb{C}}_{\sigma} ~&=\{\gamma\in \Lambda G^{\mathbb{C}}_{\sigma}~
|~ \gamma \hbox{ extends holomorphically to } |\lambda|<1 \},\\[1mm]
 \Lambda_{\mathcal{C}}^{+} G^{\mathbb{C}}_{\sigma} ~&=\{\gamma\in
\Lambda^+ G^{\mathbb{C}}_{\sigma}~|~   \gamma(0)\in  (K^\C)^0 \},\\[1mm]
\Lambda_{\mathcal{C}}^{-} G^{\mathbb{C}}_{\sigma} ~&=\{\gamma\in
\Lambda^- G^{\mathbb{C}}_{\sigma}~|~   \gamma(0)\in  (K^\C)^0 \}.\\[1mm]
\end{array}\end{equation*}

\begin{theorem}
{\em(Birkhoff and Iwasawa Decomposition Theorem for $ (\Lambda {G}_\sigma^\C )^0$)} \label{thm-birkhoff-0}
\begin{enumerate}
\item (Birkhoff Decomposition)
\begin{enumerate}
\item $ (\Lambda {G}_\sigma^\C )^0= \bigcup \Lambda^{-}_{\mathcal{C}} {G}^{\mathbb{C}}_{\sigma} \cdot \omega \cdot \Lambda^{+}_{\mathcal{C}} {G}^{\mathbb{C}}_{\sigma}$,
disjoint union,
where the $\omega$'s are representatives of the double cosets.

\item The multiplication $\Lambda_{*}^{-} {G}^{\mathbb{C}}_{\sigma}\times
\Lambda^{+}_\FC {G}^{\mathbb{C}}_{\sigma}\rightarrow
\Lambda {G}^{\mathbb{C}}_{\sigma}$ is an analytic  diffeomorphism onto the
open and dense subset $\Lambda_{*}^{-} {G}^{\mathbb{C}}_{\sigma}\cdot
\Lambda^{+}_\FC {G}^{\mathbb{C}}_{\sigma}$  of $ (\Lambda G^\C_\sigma)^0$
{\em (big Birkhoff cell)}.
\end{enumerate}
\item  (Iwasawa Decomposition)
\begin{enumerate}
\item
$(\Lambda G^{\C})_{\sigma} ^0=
\bigcup \Lambda G_{\sigma}^0\cdot \delta\cdot
\Lambda_\mathcal{C}^{+} G^{\mathbb{C}}_{\sigma},$
where the $\delta$'s are representatives of the double cosets.

 \item The multiplication $\Lambda G_{\sigma}^0 \times \Lambda_\mathcal{C}^{+} G^{\mathbb{C}}_{\sigma}\rightarrow
(\Lambda G^{\mathbb{C}}_{\sigma})^0$ is a real analytic map onto the connected open subset
$ \Lambda G_{\sigma}^0 \cdot \Lambda_\mathcal{C}^{+} G^{\mathbb{C}}_{\sigma}   = \mathcal{I}^{\mathcal{U}}_e \subset \Lambda G^{\mathbb{C}}_{\sigma}$.
\end{enumerate}
\end{enumerate}
\end {theorem}
For the loop groups related to Willmore surfaces we have the following result.
\begin{theorem}
Consider the setting $G/K = Gr_{1,3}(\mathbb{R}^{n+4}_{1}) = SO^+(1,n+3)/SO^+(1,3)\times SO(n)$.
\begin{enumerate}
\item
 There exist  two different open Iwasawa cells in the connected loop group
$(\Lambda G^{\mathbb{C}}_{\sigma})^0$,  one given by $\delta = e$ and the other one by $\delta = diag(-1,1,1,1,-1,1,1,...,1) $.

\item
There exists a closed, connected, solvable subgroup
$S \subseteq (K^\C)^0$ such that
the multiplication

$\Lambda G_{\sigma}^0 \times \Lambda^{+}_S G^{\mathbb{C}}_{\sigma}\rightarrow
(\Lambda G^{\mathbb{C}}_{\sigma})^0$ is a real analytic diffeomorphism onto the connected open subset
$ \Lambda G_{\sigma}^0 \cdot \Lambda^{+}_S G^{\mathbb{C}}_{\sigma}      \subset  \mathcal{I}^{\mathcal{U}}_e \subset(\Lambda G^{\mathbb{C}}_{\sigma})^0$.
\end{enumerate}
\end{theorem}


\subsubsection{The DPW method and potentials}

\begin{theorem}{\em(\cite{DPW})}\label{DPW}
Let $\D$ be a contractible open subset of $\C$ and $z_0 \in \D$ a base point.
\begin{enumerate}
\item
Let $f: \D \rightarrow G/K$ be a harmonic map with $f(z_0)=\ee K.$
Then the associated family  $f_{\lambda}$ of $f$ can be lifted to a map $F:\D \rightarrow \Lambda G_{\sigma}$ with $F(z_0,\lambda)= \ee$.
  Then the map $F$ takes only values in
$ \mathcal{I}^{\mathcal{U}}_e \subset \Lambda G^{\mathbb{C}}_{\sigma}$. There exists a discrete subset $\D_0\subset \D$ such that on $\D\setminus \D_0$
\[F(z,\lambda)=F_-(z,\lambda)\cdot F_+(z,\lambda),\]
where
$F_-(z,\lambda)\in\Lambda_{*}^{-} G^{\mathbb{C}}_{\sigma}
\hspace{2mm} \mbox{and} \hspace{2mm} F_+(z,\lambda)\in (\Lambda^{+} G^{\mathbb{C}}_{\sigma})^0.$
Moreover $F_-(z,\lambda)$ is meromorphic in $z \in \D$ and $F_-(z_0,\lambda)=\ee$ holds and the Maurer-Cartan form $\eta$ of $F_-$
\[\eta= F_-(z,\lambda)^{-1} \dd F_-(z,\lambda)\]
is a $\lambda^{-1}\cdot\mathfrak{p}^{\mathbb{C}}-\hbox{valued}$ meromorphic $(1,0)$-form with poles at points of $\D_0$ only.

\item Conversely, let $\eta$ be a  $\lambda^{-1}\cdot\mathfrak{p}^{\mathbb{C}}-\hbox{valued}$ meromorphic $(1,0)$-form for which the solution
to the equation
\begin{equation}
F_-(z,\lambda)^{-1} \dd F_-(z,\lambda)=\eta, \hspace{5mm} F_-(z_0,\lambda)=\ee,
\end{equation}
is meromorphic on $\D$, with  $\D_0$ as set of possible poles.
Then on the open set $\D_{\mathcal{I}} = \lbrace z \in \D; F(z,\lambda)
\in \mathcal{I}^{\mathcal{U}} \rbrace$ we
define $\tilde{F}(z,\lambda)$ via the factorization
 $\mathcal{I}^{\mathcal{U}}_e =  ( \Lambda G_{\sigma})^0 \cdot \Lambda_S^{+} G^{\mathbb{C}}_{\sigma}
\subset  \Lambda G^{\mathbb{C}}_{\sigma}$:
\begin {equation}\label{Iwa}
F_-(z,\lambda)=\tilde{F}(z,\lambda)\cdot \tilde{F}_+(z,\lambda)^{-1}.
\end{equation}
 Then one obtains an extended frame
$ \tilde{F}(z,\lambda)=F_-(z,\lambda)\cdot \tilde{F}_+(z,\lambda)$
of some harmonic map from $\D_{\mathcal{I}} \subset \D$ to $G/K$ satisfying  $\tilde{F}(z_0,\lambda)=\ee$.

The above two constructions are inverse to each other.
\end{enumerate}
\end{theorem}

The $\lambda^{-1}\cdot \mathfrak{p}^{\mathbb{C}}-\hbox{valued}$ meromorphic $(1,0)$
form $\eta$ is called the {\em normalized potential} for the harmonic
map $f$ with the point $z_0$ as the reference point (\cite{DPW}).

When permitting more powers of $\lambda$, one can obtain potentials with holomorphic coefficients, which will be called {\em holomorphic potentials}. They are not uniquely determined.
\begin{theorem} \label{pot-hol} Let $\D$ be a contractible open subset of $\C$.
Let $F(z,\lambda)$ be the frame of some harmonic map
into $G/K$. Then there exists some $V_+ \in \Lambda^{+} G^{\mathbb{C}}_{\sigma} $ such that $C(z,\lambda) =
F V_+$ is holomorphic in $z$ and in $\lambda \in \mathbb{C}^*$.
Then the Maurer-Cartan form $\eta = C^{-1} dC$ of $C$ is a holomorphic $(1,0)-$form on $\D$ and it is easy to verify that $\lambda \eta$ is holomorphic for $\lambda \in \C$.
Conversely, any harmonic map  $f: \D\rightarrow G/K$ can be derived from
such a holomorphic $(1,0)$-form $\eta$ on $\D$ by the same steps as in the previous theorem.
\end{theorem}
\subsection{Potentials of Willmore surfaces}

First we get the following description of  the potentials of strongly
conformally harmonic maps.
\begin{theorem}\cite{DoWa12} \label{normalized-potential-W}
 Let $f: \D\rightarrow SO^+(1,n+3)/SO^+(1,3)\times SO(n)$ be a strongly
conformally harmonic map with $f(0)=eK$ and
$F:\D \rightarrow (\Lambda G_{\sigma})^0$ an extended frame of $f$
such that $F(0,\lambda) = I$. Then the normalized potential of $f$ with respect to the base point $z = 0$ is of the form
\begin{equation}\label{eq-potential-W}
\eta= \lambda^{-1}\eta_{-1}\dd z,\  \hbox{ with }\ \eta_{-1}=\left(
    \begin{array}{cc}
      0 & \hat{B}_1 \\
      -\hat{B}_1^tI_{1,3} & 0 \\
    \end{array}
  \right)\dd z,\hspace{3mm} \mbox{with} \hspace{3mm} \hat{B_1}^tI_{1,3}\hat{B}_1=0,
\end{equation}
where $\hat{B}_1\dd z$ is a meromorphic $(1,0)$-form on $\D$ and $0$ is not a pole of $\hat{B}_1$.

Conversely,  any  normalized potential  satisfying  \eqref{eq-potential-W} on $\D$ induces a strongly conformally harmonic map from an open subset $0 \in \D_\mathcal{I} \subset \D$ into $SO^+(1,n+3)/SO^+(1,3)\times SO(n)$.
\end{theorem}

As applications, we obtain the following descriptions of the normalized potential of  strongly conformally harmonic maps which contains a constant lightlike vector, which allow us to figure out strongly conformally harmonic maps  producing constant maps.
\begin{theorem}\cite{Wang-Min} \label{th-potential-light}  Let $f: \mathbb{D}\rightarrow SO^+(1,n+3)/SO^+(1,3)\times SO(n)$ be a strongly conformally harmonic map which contains a constant light-like vector and assume that $f(0)=I_{n+4}\cdot K$ holds. Then the normalized potential of $f$ with reference point $p$ is of the form
 \begin{equation}\label{eq-w-minimal}
\eta=\lambda^{-1}\left(
                     \begin{array}{cc}
                       0 & \hat{B}_1 \\
                       -\hat{B}^{t}_1I_{1,3} & 0 \\
                     \end{array}
                   \right)\dd z,\ \hbox{ where }\ \hat{B}_{1}=\left(
\begin{array}{cccc}
 \hat{f}_{11} & \hat{f}_{12} & \cdots &  \hat{f}_{1n} \\
 -\hat{f}_{11} &  -\hat{f}_{12} & \cdots &  -\hat{f}_{1n} \\
 \hat{f}_{31} &\hat{f}_{32} & \cdots &  \hat{f}_{3n} \\
 i\hat{f}_{31} & i\hat{f}_{32} & \cdots &  i\hat{f}_{3n} \\
 \end{array}
\right).\end{equation}
Here all  $f_{ij}$ are meromorphic functions on $\mathbb{D}$.

Conversely, let  $\eta$ be a normalized potential of the form \eqref{eq-w-minimal}. Then we obtain a strongly conformally harmonic map $f: \mathbb{D}\rightarrow SO^+(1,n+3)/SO^+(1,3)\times SO(n)$ which contains a constant light-like vector.
Here we need to make sure that $\eta$ can be integrated globally.
\end{theorem}

\section{Potentials of minimal surfaces in space forms}

It is well-known that minimal surfaces in Riemannian space forms can be characterized by the following lemma (\cite{Helein}).
\begin{lemma} \cite{Helein}, \cite{Xia-Shen} Let $y:M\rightarrow S^{n+2}$ be a Willmore surface, with $\mathcal{F}$ as its conformal Gauss map. We say that $\mathcal{F}$ contains a constant vector $\mathbf{a}\in\R^{1,n+3}$ if for any
$p\in M$, $\mathbf{a}$ is contained in the  $4-$dimensional  Lorentzian subspace $\mathcal{F}(p)$. Then
\begin{enumerate}
\item $y$ is conformally  equivalent  to a minimal surface in $\mathbb{R}^{n+2}$ if and only if $\mathcal{F}$ contains a non-zero constant lightlike vector.
\item $y$ is  conformally equivalent to a minimal surface in some $S^{n+2}(c)$ if and only if $\mathcal{F}$ contains a non-zero constant timelike vector.
\item $y$ is conformally  equivalent to a minimal surface in $\mathbb{H}^{n+2}(c)$ if and only if $\mathcal{F}$ contains a non-zero constant spacelike vector.
\end{enumerate}
 \end{lemma}
By this lemma and Wu's formula, one can describe minimal surfaces in space forms as follows.
\begin{theorem}\label{th-mini-spaceform} \cite{Helein}, \cite{Xia-Shen}, \cite{Wang-Min}
Let $\mathcal{F}:M\rightarrow SO^+(1,n+3)/SO(1,3)\times SO(n)$ be a strongly conformally harmonic map. Let
\begin{equation}\label{eq-potential}
    \eta=\lambda^{-1}\left(
    \begin{array}{cc}
      0 & \hat{B}_1 \\
      -\hat{B}_1^tI_{1,3} & 0 \\
    \end{array}
  \right)dz, \hspace{4mm}  \hat{B}_{1}=\left(
\begin{array}{cccc}
 \mathbf{v}_1 & \mathbf{v}_2 & \ldots &  \mathbf{v}_n\\
 \end{array}
\right),\
\end{equation}
be the normalized potential of $\mathcal{F}$ with respect to some base point $z_0$. Then, up to a conjugation by some $T\in O^+(1,3)\times O(n)$,
\begin{enumerate}
\item $\mathcal{F}$ contains a constant lightlike vector, if and only if every $\mathbf{v}_j$ has the form
\begin{equation}\label{eq-potential-r}
 \mathbf{v}_j=
f_{j0}\left(
                                      \begin{array}{cccc}
                                        f_{j1} &
                                        -f_{j1} &
                                        f_{j3} &
                                        if_{j3} \\
                                      \end{array}
                                    \right)^t, \hbox{ with }  f_{jl}  \hbox{ meromorphic}.
\end{equation}

\item $\mathcal{F}$ contains a constant timelike vector, if and only if every $\mathbf{v}_j$ has the form
\begin{equation}\label{eq-potential-s}
\mathbf{v}_j=
g_{j}\left(
                                      \begin{array}{cccc}
                                        0 &
                                        2g_{0}&
                                        1-g_{0}^2 &
                                        i(1+g_{0}^2)\\
                                      \end{array}
                                    \right)^t, \hbox{ with }  g_{j},\ g_0  \hbox{ meromorphic}.
\end{equation}

\item $\mathcal{F}$ contains a constant spacelike vector, if and only if every $\mathbf{v}_j$ has the form
\begin{equation}\label{eq-potential-h}
 \mathbf{v}_j=
h_{j}\left(
  \begin{array}{cccc}
   2ih_{0} &
  0 &
 1-h_{0}^2 &
  i(1+h_{0}^2)\\
    \end{array}
  \right)^t, \hbox{ with }  h_{j},\ h_0  \hbox{ meromorphic}.
\end{equation}
\end{enumerate}
Potentials of the above form will  be called  canonical potentials for the  corresponding minimal surfaces in space forms.
\end{theorem}
\begin{remark}
    So far we obtain the normalized potentials. By similar arguments, one can derive the characterization of holomorphic potentials of minimal surfaces in space forms up to a gauge or a  dressing action.
\end{remark}

\section{On totally isotropic Willmore $2$-spheres in $S^6$}

Using the frame we discussed before, we first construct a new Willmore two-sphere in $S^6$, which is totally isotropic and non-S-Willmore, that is, it does not admit a dual Willmore surface.
This answers an open problem of Ejiri in \cite{Ejiri1988}.

For totally isotropic Willmore $2$-spheres in $S^6$, we first obtain the following geometric description.
\begin{theorem} \label{thm-iso-willmore-1}
Let $y$ be a Willmore two sphere in $S^6$ with isotropic Hopf differential, i.e., $\langle\kappa,\kappa\rangle=0$. If $y$ is not S-Willmore, then $y$ is totally isotropic (and hence full) in  $S^6$. Moreover, locally there exists an isotropic frame $\{E_1, E_2\}$ of the normal bundle $V^{\perp}_{\C}$ of $y$ such that
\begin{equation}\label{eq-normal-bundle}
\left\{\begin{split}
&\kappa,\ D_{z}\kappa,\ D_{\bar{z}}\kappa\in \hbox{Span}_{\C}\{E_1, E_2\},\\
&\langle E_i, E_j\rangle=0, \langle E_i, \bar{E}_j\rangle=2\delta_{ij}, i,j=1,2.\\
&  D_{z} E_i\in \hbox{Span}_{\C}\{E_1, E_2\},\ D_{\bar{z}} E_i\in \hbox{Span}_{\C}\{E_1, E_2\},\ i,j=1,2.\\
                 \end{split}\right.
\end{equation}
That is, the normal connection is block diagonal under the frame $\{E_1, E_2, \bar{E}_1, \bar{E}_2 \}$.
\end{theorem}

The loop group proof of the above theorem is due to the following two theorems, characterizing all totally istropic Willmore  two spheres in $S^6$ via normalized potentials.
\begin{theorem}\label{thm-iso-willmore-2}
Let $y$ be a totally isotropic Willmore two sphere in $S^6$. Then the normal bundle of $y$ satisfies the properties  \eqref{eq-normal-bundle} of Theorem \ref{thm-iso-willmore-1}. The normalized potential of $y$ is of the form
\begin{equation}\label{eq-iso-np}
\eta=\lambda^{-1}\left(
\begin{array}{cc}
0 & \hat{B}_1 \\
-\hat{B}_1^tI_{1,3} & 0 \\
\end{array}
\right)\dd z, \hbox{ with } \hat{B}_1=
                                        \left(
                                          \begin{array}{ccccc}
                                            h_{11} & ih_{11}&  h_{12} & ih_{12} \\
                                            h_{21} & ih_{21}&  h_{22} & ih_{22} \\
                                            h_{31} & ih_{31}&  h_{32} & ih_{32} \\
                                            h_{41} & ih_{41}&  h_{42} & ih_{42} \\
                                          \end{array}
\right)\ \&\ \hat{B}_1^tI_{1,3}\hat{B}_1=0.\end{equation}.
                    Here $h_{ij}$ are meromorphic functions on $S^2$.
\end{theorem}

\begin{theorem} \label{thm-iso-s6-2} Let $y$ be a Willmore surface in $S^6$ with its normalized potential being of the form \eqref{eq-iso-np}.
Then $y$ is totally isotropic in $S^6$. Moreover, locally there exists an isotropic frame $\{E_1, E_2\}$ of the normal bundle $V^{\perp}_{\C}$ of $y$ such that \eqref{eq-normal-bundle} holds.
\end{theorem}

In terms of the above theorem, we also derive an algorithm to derive  concrete totally istropic Willmore  surfaces in $S^{6}$ with \eqref{eq-normal-bundle} holding.  The first new non-S-Willmore  Willmore two-sphere is given as an example:

\begin{theorem} Set \begin{equation}\label{eq-example-np}\eta=\lambda^{-1}\left(
                      \begin{array}{cc}
                        0 & \hat{B}_1 \\
                        -\hat{B}_1^tI_{1,3} & 0 \\
                      \end{array}
                    \right)dz,\ \hbox{ with } \ \hat{B}_1=\frac{1}{2}\left(
                     \begin{array}{cccc}
                       2iz&  -2z & -i & 1 \\
                       -2iz&  2z & -i & 1 \\
                       -2 & -2i & -z & -iz  \\
                       2i & -2 & -iz & z  \\
                     \end{array}
                   \right).\end{equation}
The associated family of unbranched Willmore two-spheres $x_{\lambda}$, $\lambda\in S^1$, corresponding to $\eta$, is

\begin{equation}\label{example1}
 x_{\lambda} =\frac{1}{ \left(1+r^2+\frac{5r^4}{4}+\frac{4r^6}{9}+\frac{r^8}{36}\right)}
\left(
                          \begin{array}{c}
                            \left(1-r^2-\frac{3r^4}{4}+\frac{4r^6}{9}-\frac{r^8}{36}\right) \\
                            -i\left(z- \bar{z})(1+\frac{r^6}{9})\right) \\
                            \left(z+\bar{z})(1+\frac{r^6}{9})\right) \\
                            -i\left((\lambda^{-1}z^2-\lambda \bar{z}^2)(1-\frac{r^4}{12})\right) \\
                            \left((\lambda^{-1}z^2+\lambda \bar{z}^2)(1-\frac{r^4}{12})\right) \\
                            -i\frac{r^2}{2}(\lambda^{-1}z-\lambda \bar{z})(1+\frac{4r^2}{3}) \\
                            \frac{r^2}{2} (\lambda^{-1}z+\lambda \bar{z})(1+\frac{4r^2}{3})  \\
                          \end{array}
                        \right), \ \hbox{ with $r=|z| .$ }
\end{equation} Moreover $x_{\lambda}:S^2\rightarrow S^6$ is a Willmore immersion in $S^6$, which is full, not S-Willmore, and totally isotropic.
Note that for all $\lambda\in S^1$, $x_{\lambda}$ is isometric to each other in $S^6$.
\end{theorem}

Finally we note that
Theorem \ref{thm-iso-willmore-2} and Theorem \ref{thm-iso-s6-2} can be generalized to all totally istropic Willmore  two spheres in $S^{2m}$ via normalized potentials in a similar way   with many more computations and with the results  of \cite{Ma}.

\section{Symmetric Willmore surfaces via potentials}

Surfaces with symmetries are always of importance in geometry. It is therefore natural to see how the symmetries influence the potentials of symmetric Willmore surfaces. It turns out that there are two different cases which  lead to oriented Willmore surfaces \cite{DoWa-sym1} and non-oriented Willmore surfaces \cite{DoWa-sym2} respectively.

\subsection{The orientation-preserving case}

For the  orientation-preserving case, we can obtain a potential charaterization of Willmore surfaces with various kinds of symmetries, including the ones with a simple symmetry, or with a finite order symmetry, or with $S^1$-symmetries. We refer to \cite{Do-Ha2}, \cite{Do-Ha3}, \cite{Do-Ha5} and \cite{DoWa-sym1} for more details. In particular for the  case of $S^1$-symmetries, we refer to \cite{BuKi} for a very satisfactory treatment of  $S^1$-equivariant harmonic maps.

Moreover, it is natural to ask whether for every Riemann surface $M$ and every harmonic map from $M$ to $G/K$ there does exist some potential $\tilde{\eta} $ on  $\D$ which generates the given harmonic maps and is invariant under the action of $\pi_1 (M)$. We can show this holds for the conformal Gauss map of a Willmore surface and therefore for the Wilmore surface. Note that in this case   it is preferable to use a holomorphic potential  instead of the normalized potential.
\begin{theorem}
Let $M$ be a Riemann surface and let $y: M \rightarrow S^{n+2}$ be a Willmore immersion in $S^{n+2}.$ Then
\begin{enumerate}
    \item If $M$ is non-compact, then there exists a holomorphic potential $\eta$  on $\D$ which is invariant under the action of $\pi_1(M)$ and generates $y$.

    \item  If $M$ is compact, then there exists a meromorphic potential $\eta$  on $\D$  which is invariant under the action of $\pi_1(M)$ and generates $y$.
\end{enumerate}
\end{theorem}

\subsection{The orientation-reversing case}

In this case, let $\mu$ be an automorphism of a Riemann surface which reverses the orientation. Then we obtain a potential description of the Willmore surface under the symmetry $\mu$. We refer to Section 3 of \cite{DoWa-sym2} for more details on this. An interesting application is to derive/characterize non-oriented Willmore surfaces.  In particular, we obtain a classification of isotropic Willmore $\R P^2$ in $S^4$  by
providing all possible normalized potentials \cite{DoWa-sym2}.

Let $y:  S^2\rightarrow S^4$ be an isotropic Willmore surface.  Then the normalized potential of $y$ is of the form \cite{DoWa12}
\begin{equation} \label{eq-b1}
 \eta=\lambda^{-1}\left(
\begin{array}{cc}
 0 & \hat{B}_1\\
 -\hat{B_1}^tI_{1,3} & 0 \\
\end{array}
 \right)\dd z \hbox{ with }
\hat{B}_1=\frac{1}{2}\left(
\begin{array}{cccc}
i(f_3'-f_2')&  -(f_3'-f_2') \\
i(f_3'+f_2')&  -(f_3'+f_2')   \\
 f_4'-f_1' & i(f_4'-f_1')   \\
 i(f_4'+f_1') & -(f_4'+f_1')  \\
 \end{array}
 \right).
\end{equation}
 Here the functions $f_j,$ $1\leq j\leq4$, are meromorphic functions on $S^2$ satisfying $f_1'f_4'+f_2'f_3'=0$. The anti-holomorphic map $\mu:S^2\rightarrow S^2$ is defined by $\mu(z)=-\frac{1}{\bar{z}}$. Then we have
\begin{theorem}\label{thm-rp2} Let $y:  S^2\rightarrow S^4$ be an isotropic Willmore surface with its normalized potential as above.
\begin{enumerate}
\item Assume that $y$ has the symmetry
$y\circ\mu=y$. Then the following equations hold
 \begin{equation}\label{eq-iso-rp2}
    \left\{
    \begin{split}
&f_1(z)+(f_1(z)f_4(z)+f_2(z)f_3(z))\overline{f_4(\mu(z))}=0,\\
&f_2(z)+(f_1(z)f_4(z)+f_2(z)f_3(z))\overline{f_2(\mu(z))}=0,\\
&f_3(z)+(f_1(z)f_4(z)+f_2(z)f_3(z))\overline{f_3(\mu(z))}=0,\\
&f_4(z)+(f_1(z)f_4(z)+f_2(z)f_3(z))\overline{f_1(\mu(z))}=0.\\
\end{split}\right.
\end{equation}
\item Conversely,  consider a normalized potential $\eta$ of the form   \eqref{eq-b1} with   meromorphic functions $f_j$. Assume that $f_j$ satisfies
  \eqref{eq-iso-rp2} and that $\eta$ has a meromorphic antiderivative.
 Then the harmonic map $f|_{\lambda=1}$ induced by $\eta$ can be projected to a map $y$ from $S^2$ to $S^{n+2}$.
 Moreover, we have
 \begin{enumerate}
\item
If $y(z,\bar{z},\lambda)$ is not a constant map, then we have
\[y(\mu(z),\overline{\mu(z)},\lambda)=y(z,\bar{z},\lambda),\] and  the harmonic map
$f|_{\lambda=1}$  is the oriented conformal Gauss map of the Willmore surface $y|_{\lambda=1}$.

\item Moreover, $y$ is conformally congruent to some minimal $\R P^2$ if and only if there exists some none-zero time-like vector $\mathrm v\in\R^4_1$ such that $v^tI_{1,3}B_1=0$. In particular, $y$ is conformally congruent to some minimal $\R P^2$ if $f_2=f_3$.
\end{enumerate}
\end{enumerate}
\end{theorem}

Set
\begin{equation}\label{eq-f_ij}
  f_1=-2mz^{2m+1},\ f_2= f_3=i\sqrt{4m^2-1}z^{2m},\ f_4=-2mz^{2m-1},
\end{equation}
with $m\in \mathbb{Z}^+$.
Then the functions $f_j, j=1,2,3,4,$ satisfiy  \eqref{eq-iso-rp2}. When $m=1$, we obtain the well-known Veronese $\R P^2$. When $m=2$, we obtain a minimal  $\R P^2$ in $S^4$ with area $10\pi$ (Here $r:=|z|$):
\begin{equation}\label{eq-non--ori-example2}
y(z,\bar{z},\lambda)=\frac{1}{  (3r^4-4r^2+3)(1+r^2)^2 \\}\left(
                       \begin{array}{c}
                         -(3r^8-8r^6+8r^4-8r^2+3) \\
                         \sqrt{15}(z+\bar{z})(1-r^2)(1+r^4) \\
                        -i\sqrt{15}(z-\bar{z})(1-r^2)(1+r^4) \\
                        \sqrt{15}(\lambda^{-1}z^4+\lambda\bar{z}^4) \\
                         i\sqrt{15}(\lambda^{-1}z^4-\lambda\bar{z}^4) \\
                       \end{array}
                     \right).
\end{equation}
We refer to \cite{DoWa-sym2}, \cite{DoWa-equ} and \cite{WW} for more developments.

\section{The Bj\"{o}rling problem for Willmore surfaces}

In \cite{BW}, Brander and the second author proposed the Bj\"{o}rling problem for Willmore surfaces and solved it in a quite general setting. Here we recall a general case for the Bj\"{o}rling problem for Willmore surfaces in $S^{n+2}$.
\vspace{2mm}

{\bf Question}: {\em Given two real analytic curves $[Y_0]$ and $[\hat{Y}_0]$ such that they are exactly the intersection of two sphere congruences $\Phi_0$ and $\hat{\Phi}_0$. Does there exist a pair of Willmore surfaces $[Y]$ and $[\hat{Y}]$, adjoint to each other, such that $[Y]$ passes $[Y_0]$ and  $[\hat{Y}]$ passes $[Y_0]$, with their conformal Gauss maps congruent to the two sphere congruences $\Phi_0$ and $\hat{\Phi}_0$ along the curves respectively?}

\begin{theorem} $($Theorem 8.8 of \cite{BW}$)$

Let $$\Phi_0=Y_0\wedge\hat Y_0\wedge P_{1}\wedge P_{2},\ \hat\Phi_0=\hat Y_0\wedge Y_0\wedge \hat{P}_{1}\wedge \hat{P}_{2}:\mathbb{I}\rightarrow SO^+(1,n+3)/(SO^+(1,3)\times SO(n))$$ denote two real analytic sphere congruences from an open unit interval $\mathbb{I}$ to $S^{n+2}$ such that
  \begin{enumerate}
\item $Y_0:\mathbb{I}\rightarrow \mathcal{C}_+^{n+3}\subset\mathbb{R}^{n+4}_1$ is a real analytic curve with arc-parameter $u$ and $[Y_0]$ is an enveloping curve of $\Phi_0$;
\item  The  real analytic map $\hat{Y}_0:\mathbb{I}\rightarrow \mathcal{C}_+^{n+3}$ satisfies  $\langle Y_0, \hat{Y}_0\rangle=-1$. And it is an enveloping curve of $\hat\Phi_0$ at the points it is immersed.
 \end{enumerate}

 Then there exists a unique isotropic harmonic map  $Y\wedge \hat Y :\Sigma\rightarrow SO^+(1,n+3)/(SO^+(1,1)\times SO(n+2))$ and a unique Willmore surface $y=[Y]:\Sigma\rightarrow  S^{n+2}$, with   an adjoint transform $\hat y=[\hat Y]$, $\Sigma$ some simply connected open subset containing $\mathbb{I}$ and $z=u+iv$ a complex coordinate of $\Sigma$, such that:
  \begin{enumerate}
\item  The canonical lift $Y$ of $y$ satisfies
$Y|_{\mathbb{I}}=Y_0$;
\item The map $\hat Y$ satisfies $\hat Y|_{\mathbb{I}}=\hat Y_0$;
\item  The conformal Gauss map $Gr$ and $\widehat{Gr}$ of $y$ and $\hat y$ respectively  satisfies
\[Gr|_{\mathbb{I}}=\Phi_0,\ \ \widehat{Gr}|_{\mathbb{I}}=\hat\Phi_0.\]
\end{enumerate}\end{theorem}

The main idea of proof is an  application of the boundary potential developed by Brander and the first author \cite{BrDo}. The interesting part is that in this case,  one uses both,  the conformal Gauss map and another (local) harmonic map derived by Helein \cite{Helein}, Xia-Shen \cite{Xia-Shen} and Ma \cite{Ma2006}.  We refer to \cite{BW} for a detailed discussion.

\section{Beyond Willmore surfaces}

 The target manifold of the conformal Gauss map of a Willmore surface is a non-compact symmetric space. While global results for harmonic maps into compact symmetric spaces are difficult to generalize
directly to the case of non-compact symmetric spaces, a ``Duality Theorem"
(see just below) permits to translate at least some properties.

  In particular the properties of being of finite uniton type \cite{Uh}, \cite{BuGu}) are preserved under this duality \cite{DoWa13}. Moreover, we re-formulate the  theory of Burstall and Guest  \cite{BuGu} to algebraic harmonic maps into compact inner symmetric spaces, describing harmonic maps by potentials and fixed initial condition $I$ at a fixed base point, thus completely describing them in terms of DPW.
  This permits to  apply the Duality Theorem to obtain a characterization of harmonic $2$-spheres in non-compact inner symmetrice spaces.
  In particular, this allows us to derive all Willmore $2-$spheres in $S^{n+2}$ via normalized potentials (see \cite{Wang-1} for more details).

\subsection{ The Duality Theorem for harmonic maps}

First we state the  Duality Theorem as follows.
\begin{theorem}\label{thm-noncompact}
	Let $G/K$ be an inner, non-compact symmetric space with $G$ being semisimple and  simply connected. Let $e$ be the identity element of $G$. Then there exists a maximal compact Lie subgroup $U$ of $G^{\C}$ such that $({U} \cap  {K}^{\mathbb{C}})^{\mathbb{C}}= {K}^{\mathbb{C}}$ and $ {U}/ ( {U} \cap  {K}^{\mathbb{C}} )$ is a compact, inner  symmetric space. Moreover, let $  M$  be a connected, simply connected Riemann surface with a basepoint $p_0 \in M$.
	\begin{enumerate}
		\item Let $f:  M\rightarrow  G/ K$ be a harmonic map with    $f(p_0)=e K$.
  Then there exists a  unique harmonic map  $f_{U}: M \rightarrow   U/ (U \cap  {K}^{\mathbb{C}})$ into the compact, inner  symmetric space $ U/ ( U \cap  K^{\mathbb{C}} )$ which has the same normalized potential as $f$ and $f_U(p_0)=e ( {U} \cap  {K}^{\mathbb{C}} )$.
		\item Conversely,  let $ f_U:  {M} \rightarrow {U}/ ( {U} \cap  {K}^{\mathbb{C}} )$ be a harmonic map with $f_U(p_0)=e ( {U} \cap  {K}^{\mathbb{C}} )$.
		Then there exists a neighbourhood  $ {M}_0\subset  {M}$ of $p_0$ and a harmonic map
		$ {f} :  {M}_0 \rightarrow  {G}/ {K}$ which has the same normalized potential as
		$  {f}_{ {U}}$.
	\end{enumerate}
The harmonic maps $f$ and $f_U$ are related via the Iwasawa decompositions of their extended frames  w.r.t. different real loop groups $\Lambda G_{\sigma}$ and $\Lambda U_{\sigma}$.
\end{theorem}
We refer to  \cite{DoWa13} for more details.

As illustrations of the above theorem, we recall the following examples of harmonic maps related to Willmore surfaces in $S^4$. Note that in this case $G/K=SO^+(1,5)/SO^+(1,3)\times SO(2)$ and $U/(U\cap K^{\C})=SO(6)/SO(4)\times SO(2)$.
\begin{example} \cite{Wang-iso}
We consider  $SO^+(1,5)$ as the connected component of $SO(1,5)=\{A\in Mat(6,\R)\ | \ A^tI_{1,5}A=I_{1,5}, |A|=1\}$ containing $I$. Here $I_{1,5}=\hbox{diag}(-1,1,1,1,1,1)$. Let $f_2$ and $f_4$ be meromorphic functions on $M$ with $f_2(z_0)=f_4(z_0)=0$ at some base point $z_0\in M$, where $M$ is
a compact Riemann surface of genus $\geq 0$. Let
\begin{equation}\eta=\lambda^{-1}\left(
                      \begin{array}{cc}
                        0 &  {B}_1 \\
                        - {B}_1^tI_{1,3} & 0 \\
                      \end{array}
                    \right)dz,\ \hbox{ with }\   {B}_1=\frac{1}{2}\left(
                     \begin{array}{cccc}
                      -if_2' &  f_2'  \\
                      if_2'&  -f_2'  \\
                      f_4' & if_4'   \\
                      if_4' & -f_4'  \\
                     \end{array}
                   \right)
                   \end{equation}
                   be the normalized  potential of a  strongly conformally   harmonic map $f(z,\bar z, \lambda): M\rightarrow SO^+(1,5)/SO^+(1,3)\times SO(2)$.
                   Then the extended frame $F(z,\bar z, \lambda)=(\phi_1,\phi_2,\phi_3,\phi_4,\phi_5,\phi_6)$ of $f(z,\bar z, \lambda)$ is
                   \begin{equation}\left(
                       \begin{array}{cccccc}
                         1+\frac{|f_2|^2}{2} &\frac{|f_2|^2}{2}  & \frac{-i(\bar{f}_2f_4-f_2\bar f_4)}{2\sqrt{d_3}} & \frac{\bar{f}_2f_4+f_2\bar f_4}{2\sqrt{d_3}} &
                         \frac{-i(\lambda^{-1}f_2-\lambda\bar{f}_2)}{2\sqrt{d_3}}&   \frac{\lambda^{-1}f_2+\lambda\bar{f}_2}{2\sqrt{d_3}} \\
                         -\frac{|f_2|^2}{2} & 1-\frac{|f_2|^2}{2} & \frac{i(\bar{f}_2f_4-f_2\bar f_4)}{2\sqrt{d_3}} & -\frac{\bar{f}_2f_4+f_2\bar f_4}{2\sqrt{d_3}}  &\frac{i(\lambda^{-1}f_2-\lambda\bar{f}_2)}{2\sqrt{d_3}}&   -\frac{\lambda^{-1}f_2+\lambda\bar{f}_2}{2\sqrt{d_3}} \\
                         0 & 0&  \frac{1}{\sqrt{d_3}} & 0& \frac{\lambda^{-1}f_4+\lambda\bar{f}_4}{2\sqrt{d_3}} & \frac{i(\lambda^{-1}f_4-\lambda\bar{f}_4)}{2\sqrt{d_3}} \\
                         0 & 0 & 0 &  \frac{1}{\sqrt{d_3}}& \frac{i(\lambda^{-1}f_4-\lambda\bar{f}_4)}{2\sqrt{d_3}}&  -\frac{\lambda^{-1}f_4+\lambda\bar{f}_4}{2\sqrt{d_3}}  \\                         \frac{-i(\lambda^{-1} f_2-\lambda\bar f_2)}{2} & \frac{-i(\lambda^{-1}f_2-\lambda\bar f_2)}{2}  &-\frac{\lambda^{-1}f_4+\lambda\bar{f}_4}{2\sqrt{d_3}} & \frac{-i(\lambda^{-1}f_4-\lambda\bar{f}_4)}{2\sqrt{d_3}} &\frac{1}{ \sqrt{d_3}}&0 \\
                         \frac{\lambda^{-1}f_2+\lambda\bar f_2}{2}  & \frac{\lambda^{-1}f_2+\lambda\bar f_2}{2}    & \frac{-i(\lambda^{-1}f_4-\lambda\bar{f}_4)}{2\sqrt{d_3}} &\frac{\lambda^{-1}f_4+\lambda\bar{f}_4}{2\sqrt{d_3}} &0&\frac{1}{ \sqrt{d_3}} \\
                       \end{array}
                     \right).
                   \end{equation}
                    Here $d_1=1+|f_2|^2+|f_4|^2$, $d_2=1+|f_2|^2$ and $d_3=1+|f_4|^2$.

                    Note that $f_{\lambda}$ is the conformal Gauss map of the minimal surface $x_{\lambda}$ in $\mathbb{R}^4$:
\begin{equation}\label{eq-min-iso2}x_{\lambda}=\left(
                                             \begin{array}{cccc}
                        -\frac{if_2'}{f_4'} +\frac{i\overline{f_2'}}{\overline{f_4'}} \\
                        -\frac{f_2'}{f_4'}-\frac{\overline{f_2'}}{\overline{f_4'}} \\
                     -i(\lambda^{-1}f_2-\lambda\bar{f}_2)+\frac{i\lambda^{-1}f_2'f_4}{f_4'}-\frac{i\lambda\overline{f_2'}\bar{f}_4}{\overline{f_4'}} \\
                     (\lambda^{-1}f_2+\lambda\bar{f}_2)-\frac{\lambda^{-1}f_2'f_4}{f_4'}-\frac{\lambda\overline{f_2'}\bar{f}_4}{\overline{f_4'}} \\
                                                              \end{array}
                                           \right). \end{equation}
The extended frame  $\tilde F(z,\bar z, \lambda)$  of the dual harmonic map $\tilde f(z,\bar z, \lambda)$ into the dual compact symmetric space $SO(6)/SO(4)\times SO(2)$ is
\[\left(
                       \begin{array}{cccccc}
                         \frac{1}{\sqrt{d_2}} &0 & \frac{\bar{f}_2f_4+f_2\bar f_4}{2\sqrt{d_1d_2}} & \frac{i(\bar{f}_2f_4-f_2\bar f_4)}{2\sqrt{d_1d_2}}&   -\frac{\lambda^{-1}f_2+\lambda\bar f_2}{2\sqrt{d_1}} & -\frac{i(\lambda^{-1}f_2-\lambda\bar f_2)}{2\sqrt{d_1}}   \\
                         0 & \frac{1}{\sqrt{d_2}}& \frac{i(\bar{f}_2f_4-f_2\bar f_4)}{2\sqrt{d_1d_2}} & \frac{-\bar{f}_2f_4-f_2\bar f_4}{2\sqrt{d_1d_2}} &
                          \frac{i(\lambda^{-1}f_2-\lambda\bar f_2)}{2\sqrt{d_1}} & -\frac{\lambda^{-1}f_2+\lambda\bar f_2}{2\sqrt{d_1}}  \\
                         0 & 0&  \frac{\sqrt{d_2}}{\sqrt{d_1}} & 0 &  \frac{\lambda^{-1}f_4+\lambda\bar f_4}{2\sqrt{d_1}}    &   \frac{i(\lambda^{-1}f_4-\lambda\bar f_4)}{2\sqrt{d_1}}  \\
                         0 & 0 & 0 &   \frac{\sqrt{d_2}}{\sqrt{d_1}}  &
                          \frac{i(\lambda^{-1}f_4-\lambda\bar f_4)}{2\sqrt{d_1}} & -\frac{\lambda^{-1}f_4+\lambda\bar f_4}{2\sqrt{d_1}}  \\
                         \frac{\lambda^{-1}f_2+\lambda\bar f_2}{2\sqrt{d_2}}  &\frac{-i(\lambda^{-1}f_2-\lambda\bar f_2)}{2\sqrt{d_2}}
                         &-\frac{\lambda^{-1}f_4+\lambda\bar{f}_4}{2\sqrt{d_1d_2}} & \frac{-i(\lambda^{-1}f_4-\lambda\bar{f}_4)}{2\sqrt{d_1d_2}} & \frac{1}{\sqrt{d_1}} &0 \\
                         \frac{i(\lambda^{-1}f_2-\lambda\bar f_2)}{2\sqrt{d_2}}  &  \frac{\lambda^{-1}f_2+\lambda\bar f_2}{2\sqrt{d_2}}  & \frac{-i(\lambda^{-1}f_4-\lambda\bar{f}_4)}{2\sqrt{d_1d_2}} &\frac{\lambda^{-1}f_4+\lambda\bar{f}_4}{2\sqrt{d_1d_2}}& 0 &\frac{1}{\sqrt{d_1}} \\
                       \end{array}
                     \right).
                   \]
\end{example}

\subsection{Harmonic maps of finite uniton type}

For harmonic maps of finite uniton type, we recall the following main theorem.

\begin{theorem}
Let $\mathcal{F}: \tilde{M}\rightarrow G/K$  be a harmonic map of finite uniton type, and
 ${\mathcal{F}_ {U}:} \tilde{M} \rightarrow U/(U\cap K^{\C})$ the compact dual harmonic map of $\mathcal{F}$, with base point $z_0\in\tilde{M}$  such that $\mathcal{F}_{z=z_0}=eK$ and $(\mathcal{F}_{U})_{z=z_0}=e(U\cap K^{\C})$.  Then the normalized potential $\eta_-$ and the normalized extended framing  $F_-$  for
 $\mathcal{F}_U$,  provides also a normalized potential $\eta_-$ and a normalized extended framing $F_-$  of $\mathcal{F}$, with initial condition $F_-(z_0,\lambda)=e$. In particular, all harmonic maps of finite uniton type $\mathcal{F}: \tilde{M}\rightarrow G/K$  can be obtained in this way.
\end{theorem}

 Note that the normalized potentials for
 $\mathcal{F}_U$ can be characterized to take values in some nilpotent Lie algebra of the loop group by the theory of Burstall-Guest \cite{BuGu}. This idea has been outlined in the appendix of \cite{BuGu}. In \cite{DoWa2}, we showed that one can use the potentials stated in \cite{BuGu} to obtain all finite uniton number harmonic maps by carrying out the DPW procedure, where all frames are obtained by integrating the potentials with initial condition   $e(U\cap K^{\C})$.

\subsection{$K^{\C}$-dressing on harmonic maps and Willmore deformations}

 Dressing actions play an important role in the study of integrable system \cite{Uh}, \cite{BP}.
Usually  dressing action concerns elements in the loop group. In \cite{WW}, we realized that one can use the most simple elements, i.e., elements in $K^{\C}\setminus K$  since such elements also induce a non-trivial  dressing action on harmonic maps into $G/K$. In fact, by use of $K^{\C}$-dressing actions on the conformal Gauss map of a Willmore surface, we built an interesting Willmore deformation of minimal surfaces in $S^{n+2}$ and in $\mathbb H^{n+2}$.
\begin{theorem} $($Theorem 1.1 of \cite{WW}$)$  Let $y:U\rightarrow S^{n+2}$ be a minimal surface from a simply connected open subset $U\subset M$. There exists a family of Willmore surfaces $y_t:U'\subset U\rightarrow S^{n+2}$, $t\in[0,\pi]$, such that $y_{t}|_{t=0}=y$ and $y_{t}|_{t=\pi/2}$ is conformally equivalent to a minimal surface in  $\mathbb H^{n+2}$. Here $U'$ is an open subset of $U$.
\end{theorem}

Applying this deformation to the Veronese surface in $S^4$ and its generalization, one can derive complete minimal surfaces in $\mathbb H^4$ with any positive Willmore energy. We refer to \cite{WW} for the explicit expressions of these examples.

\section{Further developments}
\subsection{On Willmore $2$-tori in $S^n$}
The celebrated work of  Marques and Neves solved the Willmore conjecture for tori in $S^3$ and lead to
much progress in geometry and analysis \cite{Marques}. We also refer to \cite{Ku-Sch,Riviere} for related progress on the geometry and analysis of Willore surfaces.
For further understanding of Willmore tori, we still have a long road to go. It is of course an important task to consider  Willmore $2$-tori in $S^n$ via the above structure. A key difficulty is how to describe the double periodicity condition, i.e., the monodromy matrices. Neverless, related to this direction there has been much progress. For example, in \cite{Ba-Bo}, Babich and Bobenko successfully constructed Willmore tori via minimal surfaces in $H^3$.
In \cite{Pinkall1985}, Pinkall  constructed  Willmore Hopf tori by use of the Hopf bundle. Moreover, Ferus and Pedit \cite{FP} constructed more equvariant Willmore tori by use of the work of Lawson and Hsiang on minimal surfaces. Using the integrable system methods, Bohle proved that a constrained Willmore torus in $S^4$ is either of finite uniton type or of finite type in \cite{Bo}. It is an important open problem to understand his proofs in terms of the framework of \cite{DoWa12} and \cite{DoWa11}. In particular, could we derive the potentials of all $S^1$-equivariant Willmore tori in $S^3$ and $S^4$. Moreover, could we estimate their Willmore energy?

\subsection{On Willmore Klein bottles}

In \cite{Lawson}, Lawson also constructed many minimal Klein bottles in $S^3$. It is natural to ask whether there are also Willmore  Klein bottles in $S^3$ which are not minimal in any space form. A simple idea is to generalize the Lawson minimal Klein bottles in $S^3$ along the line of Ferus and Pedit's construction of  equvariant Willmore tori in $S^3$ \cite{FP}.

 Another approach could be to construct first some equvariant Willmore Mobius strip in $S^3$. This is presented  in  \cite{DoWa-sym2}, \cite{DoWa-equ}. It seems possible to us to prove the existence of Willmore Klein bottles among the  Moebius strips mentioned above.  It will be a challenging task to obtain progress in this direction.

\subsection{On higher genus cases}

For higher genus case, the study is extremely difficult. {The conjecture by Kusner \cite{Kusner} states} that the Lawson minimal surface $\xi_{1,g}$ \cite{Lawson} minimizes the Willmore energy uniquely among all surfaces in $S^n$ of genus $g$ up to a conformal transform of $S^n$.  When $g=1$, this is the (by now proven) Willmore conjecture. {So far there is very little progress on the Kusner Conjecture for $g>1$.}

It is a natural idea to use DPW to study $\xi_{1,g}$.
But it turns out to be very difficult. In \cite{HellerS1} Heller first successfully gave a DPW description for the  Lawson minimal surface $\xi_{1,2}$ in $S^3$ of genus $2$ \cite{Lawson}. Recently there has been  essential progress in this direction.   In  \cite{HHT}  Heller-Heller-Traizet obtained a finer area estimate of $\xi_{1,g}$ via the DPW potential. And in \cite{HH} Heller-Heller derived a method to conformally parametrize all Lawson surfaces $\xi_{1,g}$. There is still a long road to go for further development on  higher genus minimal surfaces or Willmore surfaces via DPW. But  the work listed above  sheds some light in  this direction and one can hope that there will be more and  more  progress along this line,  since the study of  higher genus minimal surfaces or Willmore surfaces is  extremely difficult, no matter what method will be used.

 {\bf Acknowledgements}\ \
 The second named author was partly supported by the Project 12371052 of NSFC.

\def\refname{Reference}

\vspace{2mm}
{\small
Josef Dorfmeister

Fakult\" at f\" ur Mathematik,

TU-M\" unchen, Boltzmannstr. 3,

D-85747, Garching, Germany

{\em E-mail address}: dorfm@ma.tum.de\vspace{1mm}

Peng Wang

School of Mathematics and Statistics, FJKLMAA,

Key Laboratory of Analytical Mathematics and Applications (Ministry of Education),

Fujian Normal University, Qishan Campus,

Fuzhou 350117, P. R. China

{\em E-mail address}: {pengwang@fjnu.edu.cn}
}
\end{document}